\documentclass[12pt]{article}
\usepackage{amssymb}

\newcommand{\Z}{{\mathbb Z}}

\newcommand{\eps}{{\varepsilon}}

\newtheorem{theorem}{Theorem}
\newtheorem{lemma}{Lemma}

\title{Sum-product inequalities with perturbation}

\author{Spencer Backman \\ Ernie Croot \\ Mariah Hamel \\ Derrick Hart}

\begin{document}

\maketitle

\section{Introduction}

The standard sum-product inequality, as developed by Erd\H os and Szemer\' edi \cite{erdos},
and improved upon by Elekes \cite{elekes}, Ford \cite{ford},
Nathanson \cite{nathanson}, and Solymosi \cite{solymosi1} \cite{solymosi2},
asserts that if $A$ is a set of $n$ real numbers (though Erd\H os and Szemer\' edi
originally prove their theorem for $A \subseteq \Z$), then
$$
|A+A|\ +\ |A.A|\ \geq\ n^{1 + \eps},\ {\rm whenever\ } n > n_0.
$$

It is interesting to consider whether one can prove a finer version of this theorem,
where one is allowed to ``perturb'' the products $ab \in A.A$ a little bit:
suppose that for each pair $a,b \in A$ we get to choose
``perturbation parameters'' $\delta_{a,b}$ and $\delta'_{a,b}$, where say
$$
\delta_{a,b}, \delta'_{a,b}\ \in\ [-1,1].
$$
Given these parameters, define the ``perturbed product set''
$$
P\ :=\ \{ (a + \delta_{a,b})(b + \delta'_{a,b})\ :\ a,b \in A\}.
$$
Must it be the case that for all choices of $\delta_{a,b}, \delta'_{a,b}$
we get that
$$
|A+A| + |P|\ \geq\ n^{1+\eps}\ ?
$$

The answer is obviously `no', since if the elements of $A$ are close enough together,
and are in arithmetic progression, we can easily choose $\delta_{a,b}$ and $\delta'_{a,b}$ for
every $c \in A$, so that
$$
|A+A| + |P|\ =\ 2n.
$$
But what if we add the condition that the elements of $A$ are all spaced at least $1$
apart?  Can we show that $|A+A| + |P|$ must always be large?

Again the answer is `no', but the natural example is a little more complicated.
Basically, consider the arithmetic progression
\begin{equation} \label{A_example}
A\ :=\ \{x+1, x+2, ..., x+n\}.
\end{equation}
This set has the property that it is ``close'' to the geometric progression
\begin{equation} \label{A_close}
\{ x(1 + 1/x)^j\ :\ j=0,1,...,n-1\}.
\end{equation}
Indeed,
$$
x(1 + 1/x)^j\ =\ x + j  + O(j^2/x),
$$
for $x \geq n$.  So, clearly, for $x$ big enough relative to $n$ we will have that
if the $\delta_{a,b}, \delta'_{a,b}$ are allowed to vary over $[-1,1]$, we can force
$$
P\ =\ A+A,
$$
thereby making
$$
|A+A| + |P|\ \leq\ 4n-2.
$$

But we can give even better upper bounds on $\delta_{a,b}$ and $\delta'_{a,b}$
that make $|P|+|A+A|$ small:  first, set all $\delta'_{a,b} =0$.  Then suppose that
$a$ and $b$ correspond to the numbers (\ref{A_close}) using exponents 
$j$ and $k$, respectively.  Next, we choose $\delta_{a,b} \ll n/x$ so that
\begin{eqnarray*}
(a + \delta_{a,b})b\ &=&\ x^2 (1 + j/x + O(j^2/x^2) + \delta_{a,b})(1 + k/x + O(k^2/x^2)) \\
&=&\ x^2 ( 1 + (j+k)/x + \delta_{a,b} + O(n^2/x^2)).
\end{eqnarray*}
So, if $\delta_{a,b}$ is the negative of this $O(n^2/x^2) < O(n/x)$, then
the product will be just $x^2 ( 1 + (j+k)/x)$, meaning that we can make our perturbed
product set into an arithmetic progression, making $|P| + |A+A| \ll n$.  

So, to achieve a lower bound on $|P|+|A+A|$ that is substantially better than the trivial
bound, we will need that $\delta_{a,b},\delta'_{a,b}$
can only vary over intervals that are width $O(n/x)$.  Our 
main theorem below will show that this is essentially best possible, as the intervals
for $\delta_{a,b}$ and $\delta'_{a,b}$ 
leading to non-trivial results will have lengths $n^{1-\eps}/a$ and $n^{1-\eps}/b$,
respectively.  Actually, our theorem is even more general, since following Remark 2 at the
end of our theorem, the sumset $A+A$ can also be perturbed, and still we get
a good lower bound on the resulting perturbed sum and product sets.

Our main theorem is as follows:

\begin{theorem} \label{main_theorem}  Suppose $0 < \eps \leq 1$ and let $A$
be a set of $n > n_0(\eps)$ positive real numbers, all at least $1$ apart.
For each pair $(a,b) \in A \times A$, suppose that $\delta_{a,b}$ and $\delta'_{a,b}$ 
are arbitrary real numbers satisfying
\begin{equation} \label{delta_upper}
|\delta_{a,b}|\ <\ {n^{1-\eps}  \over a},\ {\rm and\ } |\delta'_{a,b}|\ <\ 
{n^{1-\eps}  \over b}.
\end{equation}
Finally, define the {\it perturbed product set}
$$
P\ :=\ \{ (a + \delta_{a,b})(b + \delta'_{a,b})\ :\ a,b \in A\}.
$$
Then, we have that
$$
|P| + |A+A|\ \gg\ {n^{1+\eps/9} \over  \log n}.
$$
\end{theorem}

\noindent {\bf Remark 1.}  Obviously, if the elements of $A$ are not at least $1$ apart, we can rescale to
make it true.  Furthermore, all we really need is that the median of the gaps between consecutive elements of $A$
is at least $1$, since by deleting at most $n/2$ elements from $A$ we get a set of elements that are all
at least $1$ apart.
\bigskip

\noindent {\bf Remark 2.} From the proof one can show that if we also perturb the sums $A+A$,
we have the same quality lower bound on the sums and products;
that is, suppose we define $S$ to be the set
of all perturbed sums $a + b+ \delta''_{a,b}$, where $|\delta''_{a,b}| \leq n^{1-\eps}/(a+b)$.  Then, 
we can show
\begin{equation} \label{product_sum_perturb}
|P| + |S|\ \gg\ {n^{1+\eps/9} \over \log n}.
\end{equation}
Basically, the reason we can show this is that in the first parts of the proof, we pass to a subset 
$B \subseteq A$, contained in a dyadic interval $[x,2x)$, whose set of perturbed products or sums $B+B$
we show must contain at least $n^{1+\eps/9}/3\log n$ elements.  This then means that 
$|\delta''_{a,b}| < n^{1-\eps} /x$; and, then we can bound $|S|$ from below
by a constant multiple of the set of sums of $C+C$, where $C$ is the set of elements of $B$ rounded to
the nearest multiple of $n^{1-\eps}/x$.  Rewriting the perturbed products for $B$ in terms of 
perturbed products for $C$, it is easy to see that this implies (\ref{product_sum_perturb}).
\bigskip

\section{Proof of Theorem \ref{main_theorem}}

\subsection{Preliminaries}

We will basically follow a variant of Elekes's original argument used to prove
that if $A$ is a set of $n$ reals, then
$$
|A+A| \cdot |A.A|\ \gg\ n^{5/2},
$$
from which it follows that
$$
|A+A| + |A.A|\ \gg\ n^{5/4}.
$$
But our approach will differ in that the
Szemer\'edi-Trotter theorem \cite{szemeredi} is not directly amenable to our particular
approach.  Instead, we apply a very minor generalization of the Szemer\'edi-Trotter curve
theorem of Sz\'ekely \cite{szekely} (hardly any generalization at all), which follows by the same proof as that of
Sz\'ekely.  

\begin{theorem} \label{sztr}  Suppose that one has a collection of
$\ell$ non-self-crossing curves and $p$ points.  Let $C$ denote the collections of curves.
Let $m_1$ denote the maximal number of curves that
can pass through any given pair of points $p_1,p_2$, and let $m_2$ denote
the ``average intersection multiplicity'', defined as follows
$$
m_2\ :=\ {\ell \choose 2}^{-1} \sum_{ \{c_1,c_2\} \subseteq C} |c_1 \cap c_2|.
$$
Then, the number of incidences $I$, which is the number of point-curve pairs, where
the point is on the curve, satisfies
$$
I\ \ll\ (m_1 m_2)^{1/3} (p \ell)^{2/3} + \ell + m_1 p.
$$
\end{theorem}

The way this differs from the Szemer\'edi-Trotter curve theorem in \cite{szekely} is
that $m_2$ is the {\it average} intersection multiplicity among the curves, not an absolute
upper limit on the intersection multiplicity between pairs of curves.

The proof of Szekely begins with the following result on the crossing number
${\rm cr}(G)$ of a multigraph $G$, as appears in 
\cite[Theorem 7]{szekely}.

\begin{theorem}  Suppose that $G$ is a multigraph with $n$ nodes, $e$ edges
and edge multiplicity $m$.  Then, either $e < 5nm$ or ${\rm cr}(G) \geq c e^3 /(n^2 m)$.
\end{theorem}

\noindent {\bf Proof of Theorem \ref{sztr}.}
Now to prove Theorem \ref{sztr} we construct a graph as follows:  fix one of our $\ell$
curves, and consider which of our $n$ points happens to lie on it.  By choosing a direction with
which to traverse the curve, we create an ordering of these incident points.  If there are $x$
such points on the given curve, then we form $x-1$ curve segments that adjoin consecutive
pairs of points.  We throw away the ``infinite parts'' of the curves as they will play no further role
in our proof.  

Letting $I$ denote our total number of incidences, and $e$ the number of edges in our graph, 
we will have
$$
e\ =\ I - \ell,
$$
since the number of edges each curve contributes is one less than its number of incident points.
Letting $C$ denote our set of curves, we also have that
$$
{\rm cr}(G)\ \leq\ \sum_{ \{c_1,c_2\} \subset C} |c_1 \cap c_2|\ =\ {\ell \choose 2} m_2.
$$
Finally, note that in our drawing of the graph $m = m_1$.  

Putting all this together, we either have that $e < 5m_1n$, which would imply
$$
I\ =\ e + \ell\ <\ 5m_1 n + \ell,
$$
which implies our theorem, or else
$$
(I - \ell)^3/(n^2 m_1)\ =\ e^3/(n^2 m_1)\ \ll\ {\rm cr}(G)\ \leq\ {\ell \choose 2} m_2.
$$
Theorem \ref{sztr} is now proved.
\hfill $\blacksquare$

\subsection{Restricting to a dyadic interval}

Later, we will require an upper bound on the number of curves passing
through pairs of grid points, and to achieve such upper bounds it will be good to
first pass to elements of $A$ that lie in a dyadic interval.  To this end we 
will require the following lemma.

\begin{lemma} \label{dyadic_lemma}  Suppose that $A$ is a set of $n$ real numbers satisfying
$|A+A| \leq n^{1+\delta}/3\log n$.  Then, there exists a dyadic interval $[x,2x)$ containing
at least $n^{1-\delta}$ elements of $A$.
\end{lemma}

A version of this lemma can be proved without too much trouble using only
very elementary ideas; however, we give a proof using the 
Ruzsa-Plunnecke inequality, since it makes the proof short and 
transparent.  First, let us state the Ruzsa-Plunnecke inequality \cite{ruzsa}.

\begin{theorem} \label{ruzsa_plunnecke}  Suppose that $A$ is a subset of an 
additive abelian group, such that
$$
|A+A|\ \leq\ K|A|.
$$
Then,
$$
|kA - \ell A|\ =\ |A+A+\cdots + A - A - A - \cdots - A|\ \leq\ K^{k+\ell} |A|.
$$
\end{theorem}
\bigskip

\noindent {\bf Proof of Lemma \ref{dyadic_lemma}.}  Suppose, for proof by contraposition, that every
dyadic interval contains fewer than $n^{1-\delta}$ elements of $A$.  Then, it requires at
least $n^\delta$ disjoint dyadic intervals to contain all the elements of $A$, and therefore
choosing one element from every other dyadic interval (if these disjoint dyadic intervals
are put into increasing order), we get a sequence of at least $n^\delta/2$ elements
$$
A'\ :=\ \{a_1, ..., a_k\}\ \subseteq\ A,\ k > n^\delta/2,
$$
such that
$$
a_{i+1}/a_i\ \geq\ 2.
$$

It is a simple matter to show that all the $k$-fold sums of distinct elements of 
$A'$ are distinct (think about the usual proof that binary number representations are 
unique); and so, if we assume that $|A+A| = K |A|$, then from Theorem
\ref{ruzsa_plunnecke}
$$
(n^\delta/k)^k\ \leq\ {n^\delta/2 \choose k}\ \leq\ |kA'|\ \leq\ |kA|\ \leq\ K^k n.
$$
It follows that
$$
K\ \geq\ n^{\delta - 1/k}/k.
$$
Choosing $k \sim \log n$ we get that
$$
K\ >\ n^\delta/3\log n.
$$
But this means that 
$$
|A+A|\ =\ K n\ >\ n^{1+\delta}/3 \log n,
$$
so the lemma is proved.
\hfill $\blacksquare$
\bigskip

Now we choose the dyadic interval $[x,2x)$ containing the most elements of $A$.
We may assume that this interval contains at least $n^{1-\kappa \eps}$ elements of $A$,
where we define the constant
$$
\kappa\ =\ 1/9,
$$
since otherwise Lemma \ref{dyadic_lemma}, with $\delta = \kappa \eps$ implies that 
$|A+A| > n^{1+ \kappa \eps}/3 \log n$,
thereby proving Theorem \ref{main_theorem}.

Let $B$ denote those $\geq n^{1-\kappa \eps}$ elements of $A$ contained in $[x,2x)$.
We note that if we consider the products
$$
(a + \delta_{a,b})(b + \delta'_{a,b}),\ a,b \in B,
$$
then as we vary over all legal choices for $\delta_{a,b}$ and $\delta'_{a,b}$,
since $a$ and $b$ lie in the same dyadic interval $[x,2x)$, we have that the possible
values of this product must lie in an interval of width at most 
$$
2 n^{1-\eps}  + n^{2-2\eps}/x^2.
$$
Since the gap between elements of $A$ is at least $1$, we have $x \gg n$, making our interval of
width
$$
\ll\ \Delta\ :=\ n^{1-\eps}.
$$
So, the number of our ``perturbed products'' is at least a constant multiple of the 
set of distinct values $<ab>$, $a,b \in B$, where $<t>$ denotes $t$ 
rounded to the nearest multiple of $\Delta$.

\subsection{A family of polygonal curves}

In order to apply this theorem, we need to define some points and curves that are
relevant to our problem:  we begin by letting $X$ be the elements of $B+B$;
and then we let $Y$ be the elements of $B.B$ rounded
to the nearest multiple of $\Delta$.  Our set of points will then be
$X \times Y$; so, there will be $|X|\cdot |Y|$ of them in total.  

Now we define the curves:  we begin by considering, for each $a,b \in B$, the set of
points on the line
$$
y\ =\ a(x-b),
$$
with $x \in B+b$.  Then, we round the $y$ coordinate to the nearest multiple
of $\Delta$.  Sweeping from left-to-right across the
grid, we connect consecutive points by line segments.

\subsection{Perturbing the curves} \label{perturbed_curves}

At this point we could have that some pairs of curves intersect in a segment, and
therefore have infinitely many points of intersection.  But this is easily fixed by
making a small perturbation to the curves, replacing the shared segments by closely
drawn curves that are nearly parallel and only intersect at grid points.  We furthermore
can assume that if a pair of points is common to two or more curves, then those points
must be grid points (again, by perturbing the curves slightly, while still connecting the
same grid points as before).

\subsection{The average crossing number}

Now we calculate the average intersection multiplicity of pairs of curves:  first,
we observe that the polygonal arcs are at most a vertical distance $\Delta$
from their corresponding straight lines.  And therefore, two of these polygonal arcs,
corresponding, say, to the curves
$$
y\ =\ a(x-b)\ \ {\rm and\ \ } y\ =\ c(x-d),\ c > a,
$$
can only cross if $x$ is such that
$$
|c(x-d) - a(x-b)|\ <\ 2 \Delta.
$$
In other words,
$$
|(c-a) x + ab - cd|\ <\ 2 n^{1-\eps}.
$$
Furthermore, between any consecutive $x$-values of the set
\begin{equation} \label{the_set}
\{b,c\} + B,
\end{equation}
there can be at most one crossing between the pair of polygonal arcs.

Now, since in any three consecutive points of (\ref{the_set}), two must either
belong to $b+B$ or to $c+B$, we have that every other point of (\ref{the_set}) is
at least $1$ apart.  It follows, therefore, that the number of crossings
between the two polygonal arcs is at most
$$
1 + 4 n^{1-\eps}/(c-a).
$$
(The $1$ here accounts for the ``left-most point of intersection'', and once we are given
this point, there can be at most $4n^{1-\eps}/(c-a)$ other intersection points to
the right of it.)
Thus, the average intersection multiplicity, is easily seen to be bounded from above by
\begin{eqnarray*}
1 + {2 n^{1-\eps} \over n^2} \sum_{a,c \in B \atop c > a} {1 \over c-a}
\ &\leq&\ 1 + {2 \over n^{1+\eps}} \sum_{1 \leq i < j \leq n} {1 \over (j-i)} \\
&\leq&\ 1 + {2  \log n \over n^\eps}.
\end{eqnarray*}

\subsection{The number of curves meeting in a pair of points}

Next, we need to produce an upper bound on the number of curves that can meet in a pair
of points.  We begin by noting, by the way we perturbed the curves in subsection
\ref{perturbed_curves}, that in order for two or more curves to meet in a pair of points,
those points must both be elements of $B+B$; furthermore, if any of such curves corresponds,
say, to a line $y = a(x-b)$, then that pair of points must lie in $B+b$.  And so, the pair of
points must be at least $1$ apart.

It is easy to see, then, that the curves meeting in a pair of points, $(x_1,y_1), (x_2,y_2)$,
$x_2 > x_1 + 1$, are at most in number the set of lines of the form
$$
f_i(x)\ =\ a_i (x-b_i),
$$
where all the $f_i(x_1)$ come within $\Delta$ of one another, and the same should hold for
$f_i(x_2)$.  But this implies that
\begin{equation} \label{meet_bound}
|(a_i - a_j)(x_1  - x_2)|\  \leq\ 2\Delta.
\end{equation}

First, let us see that no two of these lines can have the same slope:  if they did, say
$a_i = a_j$, then just considering the contribution of the point with $x=x_1$, we would 
have that
$$
a_i(x_1 - b_i)\ =\ a_j(x_1 - b_j),
$$
and therefore $b_i = b_j$.  But this can only hold if the two lines are in fact that same,
so we may assume the slopes are different.  Assuming this, we find from (\ref{meet_bound}) 
that
$$
|a_i - a_j|\ \leq\ 2\Delta\ =\  2 n^{1-\eps}.
$$
It is clear that, since the $a_i \in B$ are all at least $1$ apart, there can be at most 
$O(n^{1-\eps})$ choices for the slopes $a_i$ such that all pairwise differences $|a_i-a_j|$
satisfy this bound.  

We have therefore proved that 
$$
m_1\ \ll\ n^{1-\eps}.
$$
Here we are assuming that $\eps \leq 1$, because of course we know that $m_1 \geq 1$.

\subsection{Conclusion of the proof}

Putting everything together, since our $|B|^2 = n^{2-2\kappa \eps}$ lines hit the grid
$X \times Y$ in $|B| = n^{1-\kappa \eps}$ points each, we have that the number of incidences is 
$n^{3-3\kappa \eps}$.  Yet, from Theorem \ref{sztr} we find that the number of incidences is
$$
\ll\ (m_1 m_2)^{1/3} (|X|\cdot |Y| \cdot n^{2-2\kappa \eps})^{2/3}\ \ll\  n^{1/3-\eps/3} 
(|X| \cdot |Y| \cdot n^{2-2\kappa \eps})^{2/3}.
$$
So, 
$$
|X|\cdot |Y|\ \gg\ n^{2 - \eps (5 \kappa/2 - 1/2)}  
$$
So, 
$$
|X| + |Y|\ \gg\ n^{1 + \eps(1/4 - 5\kappa/4)}\ \gg\ n^{1 + \eps/9}.
$$
Note that here is where we used the fact that $\eps$ is sufficiently small -- it allowed us to ignore
the contribution of the terms $\ell + m_1 p$.

This completes the proof.

\end{document}